# On the Relation of Schatten Norms and the Thompson Metric

David A. Snyder

**Abstract**

The Thompson metric provides key geometric insights in the study or non-linear matrix equations and in many optimization problems. However, knowing that an approximate solution is within $d_T$ units of the actual solution in the Thompson metric provides little insight into how good the approximation is as a matrix or vector approximation. That is, bounding the Thompson metric between an approximate and accurate solution to a problem does not provide obvious bounds either for the spectral or the Frobenius norm, both Schatten norms, of the difference between the approximation and accurate solution. This paper reports such an upper bound, namely that $\|X - Y\|_p \leq 2^{\frac{1}{p}} \frac{(e^d-1)}{e^d} \max[\|X\|_p, \|Y\|_p]$ where $\|\cdot\|_p$ denotes the Schatten p-norm and $d$ denotes the Thompson metric between $X$ and $Y$. Furthermore, a more geometric proof leads to a slightly better bound in the case of the Frobenius norm, $\|X - Y\|_2 \leq \frac{(e^d-1)}{\sqrt{e^{2d}+1}} \sqrt{\|X\|_2{}^2 + \|Y\|_2{}^2} \leq 2^{\frac{1}{2}} \frac{(e^d-1)}{\sqrt{e^{2d}+1}} \max[\|X\|_p, \|Y\|_p]$.

## 1. Introduction

The Thompson metric provides key geometric insights into the study of non-linear matrix equations. In particular, many flows, which in other metrics may not even be contractions, have well-characterized contraction rates in the Thompson metric [1]. That flows arising in many non-linear optimization and control problems are contractions in the Thompson metric [2-4] endows this metric with great utility. Applications of the Thompson metric range from proofs of the existence and uniqueness of positive definite solutions for many types of non-linear equations [5] to non-linear optimization theory [3; 6].

While the Thompson metric is convenient for solving many optimization problems involving matrices, it is often more intuitive to view matrices solving such problems within more typical geometric contexts. Knowing that the solution of a problem $X$ and its n[th] approximation $X_n$ are $d_T$ units apart in the Thompson metric provides little indication of the how close $X_n$ is to $X$. I.e. knowing that $X_n \leq \alpha X$ and $X \leq \alpha X_n$ in the Löwner ordering, where $\alpha = e^{d_T}$, does not intuitively bound $\|X - X_n\|$ for any of the usual matrix norms $\|\cdot\|$. But it is $\|X - X_n\|$ in a suitable matrix norm, not $d_T$, that provides insight as to the quality of an approximation $X_n$.

In particular, considering the matrices $X_n$ and $X$ as linear operators on Euclidean vector spaces, the spectral norm, i.e. a Schatten p-norm with $p = \infty$, of $X - X_n$ is the relevant measure how well $X_n$ approximates $X$. Considering these matrices as themselves vectors in a Euclidean space, then the relevant assessment of how well $X_n$ approximates $X$ is the Frobenius norm, i.e. a Schatten norm with $p = 2$, of $X - X_n$. Therefore, it is useful to know an upper bound for the Schatten p-norm $\|X - X_n\|_p$ given some minimal information about $X$ (e.g. its norm) as well as the Thompson metric $d = d_T(X, X_n)$. However, the relation between $d_T(X,Y)$ and $\|X - Y\|_p$ (e.g. $Y = X_n$) is not well known, if it is understood at all.

This paper seeks to fill this important gap in our understanding of the relationship between Thompson metrics and Schatten norms by providing an upper bound for the Schatten p-norm $\|X - Y\|_p$ given the Thompson metric $d = d_T(X, Y)$ as well as the Schatten p-norms of X and Y. In particular, application of Weyl's inequalities establishes that $\|X - Y\|_p \leq 2^{\frac{1}{p}} \frac{(e^d - 1)}{e^d} \max[\|X\|_p, \|Y\|_p]$. Hopefully this paper will serve as the beginning of a conversation leading to ever tighter bounds on $\|X - Y\|_p$ given $d = d_T(X, Y)$ as well as minimal information about $X$ and $Y$, such as their norms and perhaps some knowledge of their spectra of eigenvalues.

## 2. Notation

This paper will generally use a consistent set of letters and symbols to denote certain matrices and their norms and eigenvalues. Let $X$ and $Y$ each denote Hermitian matrices with eigenvalues $\chi_1 \geq \ldots \geq \chi_n$ and $\upsilon_1 \geq \ldots \geq \upsilon_n$, respectively. Denote the eigenvalues of the matrix $\Delta = X - Y$ by $\delta_1 \geq \ldots \geq \delta_n$ and those of $E = -\Delta = Y - X$ by $\varepsilon_1 \geq \ldots \geq \varepsilon_n$. Note that $\delta_i = -\varepsilon_{n-i+1}$. $\|M\|$ denotes a Schatten norm of the matrix $M$ and $\|M\|_p$ specifically denotes the Schatten p - norm. Note that $\|M\|_p$ is a function of the eigenvalues $\mu_1 \geq \ldots \geq \mu_n$ of $M$: $\|M\|_p = f_p(\mu_1, \ldots, \mu_n) = (\sum_{i=1}^n |\mu_i|^p)^{1/p}$. Similarly, this paper will use the notation of $f(\mu_1, \ldots, \mu_n)$ as the functional form of $\|M\|$. Depending on the context, $\leq$ and $\geq$ denote either the usual ordering on real numbers or the Löwner ordering on matrices: i.e. $X \leq Y$ indicates that $Y - X$ is positive semidefinite. As is standard, $tr(M)$ denotes the trace of the matrix $M$.

## 3. Results

### 3.A. Proof of General Case

The proof begins with a lemma applying Weyl's inequalities to bound the eigenvalues of $\Delta = X - Y$ by the eigenvalues of $Y$ given upper and lower bounds for $X$ in the Löwner ordering. The second lemma, a consequence of the first lemma, bounds the eigenvalues of $\Delta = X - Y$ by the eigenvalues of $X$.

*Lemma 1*: Consider Hermitian matrices $X$ and $Y$ such that $X \leq \alpha Y$ and $X \geq \beta Y$. Then (A) $(\alpha - 1) \cdot \upsilon_i \geq \delta_i \geq (\beta - 1) \cdot \upsilon_i$ and (B) $|\delta_i| \leq \max[|\beta - 1|, |\alpha - 1|] \cdot \upsilon_i$.

*Proof*:

Note that $\Delta = X - \alpha Y + (\alpha - 1) \cdot Y = X - \beta Y + (\beta - 1) \cdot Y$. Let $\alpha_1$ be the maximum eigenvalue of $X - \alpha Y$ (which is negative semi-definite as $\alpha Y - X$ is positive semidefinite by the definition of $X \leq \alpha Y$) and $\beta_n$ be the minumum eigenvalue of $X - \beta Y$ (which is positive semidefinite by the definition of $X \geq \beta Y$). By hypothesis, $\alpha_1 \leq 0$ and $\beta_n \geq 0$. Note that the eigenvalues for $(\alpha - 1) \cdot Y$ and $(\beta - 1) \cdot Y$ are $|\alpha - 1| \cdot \upsilon_1, \ldots, |\alpha - 1| \cdot \upsilon_n$ and $|\beta - 1| \cdot \upsilon_1, \ldots, |\beta - 1| \cdot \upsilon_n$, respectively. By Weyl's inequalities, we have $\alpha_1 + (\alpha - 1) \cdot \upsilon_i \geq \delta_i \geq \beta_n + (\beta - 1) \cdot \upsilon_i$. Since $\alpha_1 \leq 0$ and $\beta_n \geq 0$, we have $(\alpha - 1) \cdot \upsilon_i \geq \delta_i \geq (\beta - 1) \cdot \upsilon_i$ and hence $|\delta_i| \leq \max[|\beta - 1|, |\alpha - 1|] \cdot \upsilon_i$.

*Lemma 2*: Again, consider Hermitian matrices $X$ and $Y$ such that $X \leq \alpha Y$ and $X \geq \beta Y$. Then (A) $(\frac{1}{\beta} - 1) \cdot \chi_i \geq \delta_i \geq (\frac{1}{\alpha} - 1) \cdot \chi_i$ and (B) $\delta_i \leq \max\left[\left(\frac{1}{\alpha} - 1\right), \left(\frac{1}{\beta} - 1\right)\right] \cdot \chi_i$.

*Proof*:

$X \leq \alpha Y$ and $X \geq \beta Y$ respectively imply $\frac{1}{\alpha} X \leq Y$ and $\frac{1}{\beta} X \geq Y$. Apply Lemma 1 to $Y$ (in place of $X$), $X$ (in place of $Y$), $\frac{1}{\alpha}$ (in place of $\beta$) and $\frac{1}{\beta}$ (in place of $\alpha$).

*Theorem 1*: Consider Hermitian matrices $X$ and $Y$ such that $X \leq \alpha Y$ and $X \geq \beta Y$. $\|X - Y\| \leq \min\{\max[|\alpha - 1|, |\beta - 1|] \cdot \|Y\|, \max[\left|\frac{1}{\alpha} - 1\right|, \left|\frac{1}{\beta} - 1\right|] \cdot \|X\|\}$

*Proof*:

Consider two sets of eigenvalues, $\lambda_1 \geq \cdots \geq \lambda_n$ and $\mu_1 \geq \cdots \geq \mu_n$ such that $\lambda_i \leq \mu_i$. Since $f$, a functional form of a Schatten norm, is monotonic in each variable, $\lambda_i \leq \mu_i$ implies $f(\lambda_1, \ldots, \lambda_i, \ldots, \lambda_n) \leq f(\mu_1, \ldots, \mu_i, \ldots, \mu_n)$. Given that implication and given that $f(\gamma\mu_1, \ldots, \gamma\mu_i, \ldots, \gamma\mu_n) = \gamma f(\mu_1, \ldots, \mu_i, \ldots, \mu_n)$, $|\delta_i| \leq \max[|\beta - 1|, |\alpha - 1|] \cdot v_i$, which is given by part (B) of Lemma 1, implies that $f(\delta_1, \ldots, \delta_n) \leq \max[|\alpha - 1|, |\beta - 1|] \cdot f(v_1, \ldots, v_n)$. Similarly part (B) of Lemma 2 yields $|\delta_i| \leq \max[\left|\frac{1}{\beta} - 1\right|, \left|\frac{1}{\alpha} - 1\right|] \cdot \chi_i$, which implies that $f(\delta_1, \ldots, \delta_n) \leq \max[\left|\frac{1}{\beta} - 1\right|, \left|\frac{1}{\alpha} - 1\right|] \cdot f(\chi_1, \ldots, \chi_n)$. Combining these two inequalities for $f(\delta_1, \ldots, \delta_n)$ with the definition of $f(\mu_1, \ldots, \mu_n) = \|M\|$ yields $\|X - Y\| \leq \min\{\max[|\alpha - 1|, |\beta - 1|] \cdot \|Y\|, \max[\left|\frac{1}{\alpha} - 1\right|, \left|\frac{1}{\beta} - 1\right|] \cdot \|X\|\}$

*Theorem 2*: Consider Hermitian matrices $X$ and $Y$ such that Thompson metric $d = d_T(X, Y)$ is finite. Let $\lambda_1, \ldots, \lambda_n$ denote a collection of numbers such that $|\lambda_i| = \max\{\min[|\chi_i|, |v_i|], \min[|\chi_{n-i+1}|, |v_{n-i+1}|]\} \cdot \frac{e^d - 1}{e^d}$. Then $\|X - Y\| \leq f(\lambda_1, \ldots, \lambda_n)$

*Proof*:

Let $\alpha = e^d$. By definition of the Thompson metric, we have (i) $X \leq \alpha Y$ and (ii) $Y \leq \alpha X$. Let $\beta = \alpha^{-1} = e^{-d}$. Then we have (iii) $\beta X \leq Y$ and (iv) $\beta Y \leq X$. Applying Lemma 1.A to (i) and (iii) yields $(\alpha - 1) \cdot v_i \geq \delta_i \geq (\beta - 1) \cdot v_i$. Noting that $\beta = \frac{1}{\alpha}$ and $\alpha = \frac{1}{\beta}$, applying Lemma 2.A to (i) and (iii) yields $(\alpha - 1) \cdot \chi_i \geq \delta_i \geq (\beta - 1) \cdot \chi_i$. Similarly, recall that $E = Y - X$, so reversing the roles of $X$ and $Y$ when applying Lemma 1.A, and respectively Lemma 2.A., to (ii) and (iv) yields $(\alpha - 1) \cdot \chi_i \geq \varepsilon_i \geq (\beta - 1) \cdot \chi_i$ and yields $(\alpha - 1) \cdot v_i \geq \varepsilon_i \geq (\beta - 1) \cdot v_i$, respectively. Since $\delta_i = -\varepsilon_{n-i+1}$, application of Lemmas 1 and 2 (part A) to (ii) and (iv) yield $-(\alpha - 1) \cdot \chi_{n-i+1} \leq \delta_i \leq -(\beta - 1) \cdot \chi_{n-i+1}$ and $-(\alpha - 1) \cdot v_{n-i+1} \leq \delta_i \leq -(\beta - 1) \cdot v_{n-i+1}$. Since (by the definition of the Thompson metric) $\alpha \geq 1$ and hence $(\alpha - 1) \geq (1 - \beta) = \frac{\alpha - 1}{\alpha}$, $(\beta - 1) \cdot \chi_i \leq \delta_i \leq (1 - \beta) \cdot \chi_{n-i+1}$ and $(\beta - 1) \cdot v_i \leq \delta_i \leq (1 - \beta) \cdot v_{n-i+1}$. Substituting $(1 - \beta) = \frac{\alpha - 1}{\alpha} = \frac{e^d - 1}{e^d}$. Thus, $|\delta_i| \leq \max\{\min[|\chi_i|, |v_i|], \min[|\chi_{n-i+1}|, |v_{n-i+1}|]\} \cdot \frac{e^d - 1}{e^d} = |\lambda_i|$. As in the proof of Theorem 1, the monotonicity of $f$ and definition of $|\delta_i|$ yield $\|X - Y\| \leq f(\lambda_1, \ldots, \lambda_n)$.

*Theorem 3* : Consider Hermitian matrices $X$ and $Y$ such that Thompson metric $d = d_T(X,Y)$ is finite.. Then $\|X - Y\|_p \leq 2^{\frac{1}{p}} \frac{(e^d-1)}{e^d} \max[\|X\|_p, \|Y\|_p]$.

*Proof* :

Note that $\max\{\min[|\chi_i|, |v_i|], \min[|\chi_{n-i+1}|, |v_{n-i+1}|]\} \leq \max\{\max[|\chi_i|, |v_i|], \max[|\chi_{n-i+1}|, |v_{n-i+1}|]\} = \max[|\chi_i|, |v_i|, |\chi_{n-i+1}|, |v_{n-i+1}|] = \max\{\max[|\chi_i|, |\chi_{n-i+1}|], \max[|v_i|, |v_{n-i+1}|]\}$. Since raising positive numbers to powers $\geq 1$ is monotonically increasing, a consequence of Theorem 2 is that $|\delta_i|^p \leq \frac{(e^d-1)}{e^d} \cdot \max\{\max[|\chi_i|^p, |\chi_{n-i+1}|^p], \max[|v_i|^p, |v_{n-i+1}|^p]\}$, which in turn is $\leq \frac{(e^d-1)}{e^d} \cdot \max\{|\chi_i|^p + |\chi_{n-i+1}|^p, |v_i|^p + |v_{n-i+1}|^p\}$. Hence, by the definition of and monotonicity of $f_p$, we have $\|\Delta\|_p^p \leq \left[\frac{(e^d-1)}{e^d}\right]^p \cdot \max\{\sum_{i=1}^n |\chi_i|^p + \sum_{i=1}^n |\chi_{n-i+1}|^p, \sum_{i=1}^n |v_i|^p + \sum_{i=1}^n |v_{n-i+1}|^p\}$. Since $\sum_{i=1}^n |\chi_i|^p = \sum_{i=1}^n |\chi_{n-i+1}|^p$ and $\sum_{i=1}^n |v_i|^p = \sum_{i=1}^n |v_{n-i+1}|^p$, $\|\Delta\|_p^p \leq \left[\frac{(e^d-1)}{e^d}\right]^p \cdot \max\{2\sum_{i=1}^n |\chi_i|^p, 2\sum_{i=1}^n |v_i|^p\}$, which by definition of the Schatten p-norm yields $\|X - Y\|_p^p \leq 2 \cdot \left[\frac{(e^d-1)}{e^d}\right]^p \max[\|X\|_p^p, \|Y\|_p^p]$. Taking the p'th root of both sides of the inequality yields the result.

### 3.B. The Frobenius (p = 2) Case

We begin by noting that $tr(A^T B)$ defines an inner product yielding the Frobenius norm, i.e. $\|A\|_2 = tr(A^T A)$. This, together with the commutative property of the trace, leads to the following version of the law of cosines for matrices: $\|A - B\|_2^2 = \|A\|_2^2 + \|B\|_2^2 - 2 \cdot tr(A^T B)$. Since for two (symmetric) positive semidefinite matrices $X$ and $Y$, $tr(X^T Y) = tr(XY) \geq 0$ [7], $\theta = \cos^{-1}\left[\frac{tr(XY)}{\|X\|_2 \|Y\|_2}\right] \leq 90°$ and hence $\|X - Y\|_2^2 \leq \|X\|_2^2 + \|Y\|_2^2$. Note that the Frobenius norm of a matrix is the same as the Euclidean norm of that matrix reshaped as a vector, so matrices under the Frobenius norm can be treated just as vectors in a Euclidean space.

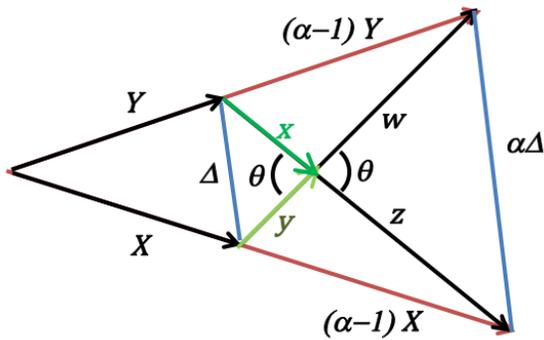

Let $d = d_T(X,Y)$ be the Thompson metric between $X$ and $Y$ and let $\alpha = e^d$. Note that in the figure at left, in which each matrix is represented as a vector, all the vectors shown are coplanar. The angle $\theta$ is the same as the angle between $\alpha X - Y = x + z$ and $\alpha Y - X = y + w$. Since, by definition of the Thompson metric, $\alpha X - Y$ and $\alpha Y - X$ are both positive semidefinite, $\theta \leq 90°$ and $tr(xy)$ and $tr(wz) \geq 0$ whereas $tr(xw)$ and $tr(yz) \leq 0$. Thus we have

(3.B.1) $$\|\Delta\|_2^2 \leq \|x\|_2^2 + \|y\|_2^2$$

(3.B.2) $$\|\alpha\Delta\|_2^2 \leq \|w\|_2^2 + \|z\|_2^2$$

(3.B.3) $$\|(\alpha - 1)X\|_2^2 \geq \|y\|_2^2 + \|z\|_2^2$$

(3.B.4) $$\|(\alpha - 1)Y\|_2^2 \geq \|w\|_2^2 + \|x\|_2^2$$

Adding Equations (3.B.1) and (3.B.2) as well as (3.B.3) and (3.B.4), we have

(3.B.5) $$\|\Delta\|_2^2 + \|\alpha\Delta\|_2^2 \leq \|x\|_2^2 + \|y\|_2^2 + \|w\|_2^2 + \|z\|_2^2 \leq \|(\alpha - 1)X\|_2^2 + \|(\alpha - 1)Y\|_2^2$$

Thus

(3.B.6) $$\|\Delta\|_2^2 + \|\alpha\Delta\|_2^2 \leq \|(\alpha - 1)X\|_2^2 + \|(\alpha - 1)Y\|_2^2$$

Solving for $\|\Delta\|_2$, we have our result $\|\Delta\|_2 \leq \frac{\alpha - 1}{\sqrt{1+\alpha^2}}\sqrt{\|X\|_2^2 + \|Y\|_2^2}$

## 4. Discussion

Weyl's inequalities, and hence some knowledge of the spectra of $X$ and $Y$, form the backbone of the proofs presented above. In the motivating case where $Y$ is an approximation of an unknown $X$, the spectrum of $X$ may also be unknown. While the principle result of this paper ultimately only requires knowledge of $\|X\|_p$ (as well as $\|Y\|_p$, which is generally known), purely geometric/trigonometric proofs, such as the one given for the Frobenius case, of the results presented in this paper would be more elegant given the nature of the motivating problem.

Furthermore, proofs not based on the matrix structure of $X$ and $Y$ but based purely on the ordering (Löwner ordering in this case) and norm (Schatten p-norm) being compared might allow for tighter bounds on $\|X - Y\|_p$ even in the absence of any knowledge of the spectrum of $X$ (or even of $Y$, for that matter), other than perhaps a restriction that $X$ and $Y$ be positive semidefinite. In comparison, Theorem 2 provides a tighter bound on $\|X - Y\|_p$ than the main result (Theorem 3), but it requires some knowledge of the spectrum of $X$ (at least that its eigenvalues are lower in magnitude than the corresponding eigenvalues of $Y$).

Additionally, proofs not based on the matrix structure of $X$ and $Y$ may lead to the generalization of these results to other orderings, which can also induce Thompson metrics [8], and other norms. For instance, since the Frobenius norm arises from an inner product, a geometrically flavored argument leads to a slightly tighter bound on $\|X - Y\|_2$ than obtained from the general bound for $\|X - Y\|_p$ and setting $p = 2$. The results presented herein may be applicable to vector spaces other than a space of square matrices. Hopefully this paper will lead to greater clarity concerning what is currently known about the relation between $\|X - Y\|_p$ and $d_T(X, Y)$ as well as spark further exploration of the relationship between metrics based on partial orders, such as the Thompson metric, and metrics induced by norms on Banach spaces.

**Acknowledgements**

This work was made possible by the support of William Paterson University of NJ's Office of the Provost for Assigned Release Time for research.